\documentclass[leqno,12pt]{amsart}
\usepackage{amssymb} 
\theoremstyle{plain} 
\newtheorem{theorem}{\indent\sc Theorem}[section] 
\newtheorem{lemma}[theorem]{\indent\sc Lemma}

\newtheorem{proposition}[theorem]{\indent\sc Proposition}

\theoremstyle{definition} 
\newtheorem{definition}[theorem]{\indent\sc Definition}

\begin{document}

\title[Double covering of the Painlev\'e I equation]{Double covering of the Painlev\'e I equation and its singular analysis \\}
\author{Yusuke Sasano }

\renewcommand{\thefootnote}{\fnsymbol{footnote}}
\footnote[0]{2000\textit{ Mathematics Subjet Classification}.
34M55; 34M45; 58F05; 32S65.}


\begin{abstract}
In this note, we will do analysis of accessible singular points for a polynomial Hamiltonian system obtained by taking a double covering of the Painlev\'e I equation. We will show that this system passes the Painlev\'e $\alpha$-test for all accessible singular points $P_i \ (i=1,2,3)$. We note its holomorphy condition of the first Painlev\'e system.
\end{abstract}
\maketitle

\section{Introduction}

In this note, we consider the first Painlev\'e equation:
\begin{equation}\label{eq:1}
{\rm P_{I}}:\frac{d^2w}{dt^2} =6w^2+t.
\end{equation}
This equation admits the following formal series expansions (see \cite{Shi}):
\begin{equation}\label{eq:sol1}
w(t)=\frac{1}{(t-t_0)^2}-\frac{t_0}{10}(t-t_0)^2-\frac{1}{6}(t-t_0)^3+h(t-t_0)^4+\frac{t_0^2}{300}(t-t_0)^6+\cdots,
\end{equation}
where $h \in {\mathbb C}$ is a free parameter and $t_0 \in {\mathbb C}$.

It is well-known that the first Painlev\'e equation is equivalent to the following Hamiltonian system (see \cite{Oka4,O3,oka6,oka7}), that is, the birational transformations
\begin{equation}\label{eq:2}
x:=w, \quad y:=\frac{dw}{dt}
\end{equation}
take the system \eqref{eq:1} into the Hamiltonian system
\begin{equation}\label{eq:3}
\frac{dx}{dt} =\frac{\partial H_{I}}{\partial y}=y, \quad \frac{dy}{dt} =-\frac{\partial H_{I}}{\partial x}=6x^2+t
\end{equation}
with the polynomial Hamiltonian:
\begin{align}\label{eq:4}
\begin{split}
&H_{I}(x,y,t)=\frac{1}{2}y^2-2x^3-tx.
\end{split}
\end{align}

It is well-known that the algebraic transformations (see P.229 \cite{1}, cf. P7 \cite{Shi}, \cite{iwasaki})
\begin{equation}\label{eq:5}
  \left\{
  \begin{aligned}
   x =&\frac{1}{v^2},\\
   y =&-\frac{2}{v^3}-\frac{1}{2}tv-\frac{1}{2}v^2+uv^3
   \end{aligned}
  \right. 
\end{equation}
take the Hamiltonian system \eqref{eq:3} to the following Hamiltonian system
\begin{equation}\label{eq:6}
  \left\{
  \begin{aligned}
   \frac{dv}{dt} &=\frac{\partial K}{\partial u}=1+\frac{t}{4}v^4+\frac{1}{4}v^5-\frac{1}{2}v^6 u,\\
   \frac{du}{dt} &=-\frac{\partial K}{\partial v}=\frac{1}{8}t^2 v+\frac{3}{8}tv^2-\left(tu-\frac{1}{4} \right)v^3-\frac{5}{4} v^4 u+\frac{3}{2} v^5 u^2
   \end{aligned}
  \right. 
\end{equation}
with the polynomial Hamiltonian
\begin{align}\label{eq:7}
\begin{split}
&K=-\frac{v^6 u^2}{4} + \frac{v^5 u}{4} + \frac{1}{4} t v^4 u - \frac{t v^3}{8} - \frac{v^4}{16} - 
 \frac{1}{16} t^2 v^2 + u.
\end{split}
\end{align}
We note that the symplectic 2-form $dy \wedge dx$ is transformed into
\begin{equation*}
dy \wedge dx=2 dv \wedge du.
\end{equation*}

Here, let us start to rewrite the Hamiltonian system  \eqref{eq:6} to a simple second-order ordinary differential equation.

The birational transformations
\begin{equation}
   q =v, \quad p =\frac{t v^4 + v^5 - 4 \frac{dv}{dt}+4}{2 v^6}
\end{equation}
take the system \eqref{eq:6} to the system
\begin{equation}
   \frac{dq}{dt} =p, \quad \frac{dp}{dt} =\frac{3}{q} p^2 - \frac{t q^3}{2} - \frac{3}{q}.
\end{equation}
We see that
\begin{equation}\label{eq:DPI}
\boxed{%
   \frac{d^2 q}{dt^2} =\frac{3}{q} \left(\frac{dq}{dt} \right)^2 - \frac{t}{2} q^3 - \frac{3}{q}.
}%
\end{equation}
For the equation \eqref{eq:DPI}, we will do the following Painlev\'e test.

Let us consider the following formal series expansions:
\begin{equation}\label{series0}
  \left\{
  \begin{aligned}
   &q=a_{-n}(t-t_0)^{-n}+a_{-(n-1)}(t-t_0)^{-(n-1)}+ \cdots+a_{-1}(t-t_0)^{-1}+a_0+ \cdots \quad (a_{-n} \not= 0, \ n \in {\mathbb N}),\\
   &\frac{dq}{dt}=-n a_{-n}(t-t_0)^{-n-1}+ \cdots,\\
   &\frac{d^2q}{dt^2}=n(n+1) a_{-n}(t-t_0)^{-n-2}+ \cdots \quad (t_0 \in {\mathbb C}).
   \end{aligned}
  \right. 
\end{equation}
Substituting the series \eqref{series0} into the equation \eqref{eq:DPI} and comparing its lowest degree, we see that
\begin{equation}\label{series10}
n=1, \quad a_{-1}=-\frac{\sqrt{2t_0}}{t_0},\frac{\sqrt{2t_0}}{t_0}.
\end{equation}
Under the conditions \eqref{series10}, we can determine its coefficients (its leading term; cf. \cite{Takei}):
\begin{equation}\label{negaseries}
q(t)=\frac{\frac{\sqrt{2t_0}}{t_0}}{t-t_0}-\frac{1}{3\sqrt{2t_0}t_0}+\frac{1}{12\sqrt{2t_0} t_0^2}(t-t_0)-\frac{5}{216\sqrt{2t_0} t_0^3}(t-t_0)^2+\frac{175 \sqrt{2} - 2592 \sqrt{2} t_0^5}{51840 \sqrt{t_0} t_0^4}(t-t_0)^3+\cdots.
\end{equation}
Making a change of variables $x:=\frac{1}{q^2} \left(=\frac{1}{v^2} \right)$ (see \eqref{eq:5}), we can obtain
\begin{equation}\label{negaseries3}
x(t)=\frac{t_0}{2}(t-t_0)^2+\frac{1}{6}(t-t_0)^3+\cdots.
\end{equation}
It is still an open question whether the above series are related to well-known 0-parameter family of formal meromorphic solutions in WKB analysis (its leading term; cf. \cite{Takei}).

On the other hand, the transformation $Q:=\frac{1}{q}$ takes the system \eqref{eq:DPI} into the equation:
\begin{equation}\label{eq:DPIB}
\boxed{%
   \frac{d^2 Q}{dt^2} =-\frac{1}{Q} \left(\frac{dQ}{dt} \right)^2 +3Q^3 + \frac{t}{2Q}.
}%
\end{equation}

For the equation \eqref{eq:DPIB}, let us consider the following formal series expansions:
\begin{equation}\label{series}
  \left\{
  \begin{aligned}
   &Q=a_{-n}(t-t_0)^{-n}+a_{-(n-1)}(t-t_0)^{-(n-1)}+ \cdots+a_{-1}(t-t_0)^{-1}+a_0+ \cdots \quad (a_{-n} \not= 0, \ n \in {\mathbb N}),\\
   &\frac{dQ}{dt}=-n a_{-n}(t-t_0)^{-n-1}+ \cdots,\\
   &\frac{d^2Q}{dt^2}=n(n+1) a_{-n}(t-t_0)^{-n-2}+ \cdots \quad (t_0 \in {\mathbb C}).
   \end{aligned}
  \right. 
\end{equation}
Substituting the series \eqref{series} into the equation \eqref{eq:DPIB} and comparing its lowest degree, we see that
\begin{equation}\label{series1}
n=1, \quad a_{-1}=-1,1.
\end{equation}
Under the conditions \eqref{series1}, we determine its coefficients:
\begin{equation}\label{formal.sol2}
  \left\{
  \begin{aligned}
   &Q(t)=\frac{1}{t-t_0}-\frac{t_0}{20}(t-t_0)^3-\frac{1}{12}(t-t_0)^4+a_5(t-t_0)^5+ \cdots \quad (a_5 \in {\mathbb C}),\\
   &Q(t)=\frac{-1}{t-t_0}+\frac{t_0}{20}(t-t_0)^3+\frac{1}{12}(t-t_0)^4+b_5(t-t_0)^5+ \cdots \quad (b_5 \in {\mathbb C}),
   \end{aligned}
  \right. 
\end{equation}
where $a_5$ and $b_5$ are free parameters.

Thus, we see that this differential equation passes the Painlev\'e test.

The system \eqref{eq:6} with \eqref{eq:7} admits a 1-parameter family of formal meromorphic solutions;
\begin{equation}\label{formal.soll}
  \left\{
  \begin{aligned}
   &v=-(t-t_0)-\frac{t_0}{20}(t-t_0)^5-\frac{1}{12}(t-t_0)^6+ \cdots,\\
   &u=\frac{4}{(t-t_0)^6}-\frac{t_0}{5(t-t_0)^2}-\frac{1}{t-t_0}+h+ \cdots,
   \end{aligned}
  \right. 
\end{equation}
where $h$ is its free parameter.

In the coordinate system $(X_1,Y_1)=(v,u v^6-4)$ (see \eqref{eq:25},\eqref{eq:26}), these solutions \eqref{formal.soll} can be rewritten as follow;
\begin{equation}
  \left\{
  \begin{aligned}
   &X_1=-(t-t_0)-\frac{t_0}{20}(t-t_0)^5-\frac{1}{12}(t-t_0)^6+ \cdots,\\
   &Y_1=t_0 (t-t_0)^4+(t-t_0)^5+h(t-t_0)^6+ \cdots
   \end{aligned}
  \right. 
\end{equation}
and in the coordinate system $(X_2,Y_2)=\left(v,u-\frac{t}{v^2}-\frac{4}{v^6} \right)$ (see \eqref{holo2}), these solutions \eqref{formal.soll} can be rewritten as follow;
\begin{equation}
  \left\{
  \begin{aligned}
   &X_2=-(t-t_0)-\frac{t_0}{20}(t-t_0)^5-\frac{1}{12}(t-t_0)^6+ \cdots,\\
   &Y_2=h+{\mathcal O}((t-t_0)),
   \end{aligned}
  \right. 
\end{equation}
where the symbol ${\mathcal O}$ denotes Landau symbol.

In the coordinate system $(z_1,w_1)=\left(\frac{1}{v},(u v-1/2)v-t/2)v^2 \right)$ (see \eqref{eq:10}), these solutions \eqref{formal.soll} can be rewritten as follow;
\begin{equation}
  \left\{
  \begin{aligned}
   &z_1=-\frac{1}{t-t_0}+\frac{t_0}{20}(t-t_0)^3+\frac{1}{12}(t-t_0)^4+h (t-t_0)^5 \cdots,\\
   &w_1=\frac{4}{(t-t_0)^2}+\frac{t_0}{10}(t-t_0)^2+\frac{1}{3}(t-t_0)^3+h(t-t_0)^4+ \cdots,
   \end{aligned}
  \right. 
\end{equation}
where this solution in $z_1$ is equivalent to the second case in \eqref{formal.sol2}.

Finally, in the coordinate system $(x,y)=\left(\frac{1}{v^2},-\frac{2}{v^3}-\frac{1}{2}t v-\frac{1}{2}v^2+u v^3 \right)$, these solutions \eqref{formal.soll} can be rewritten as follow;
\begin{equation}
  \left\{
  \begin{aligned}
   &x=\frac{1}{(t-t_0)^2}-\frac{t_0}{10}(t-t_0)^2-\frac{1}{6}(t-t_0)^3+h(t-t_0)^4+ \cdots,\\
   &y=-\frac{2}{(t-t_0)^3}-\frac{t_0}{5}(t-t_0)-\frac{1}{2}(t-t_0)^2+4 h(t-t_0)^3+ \cdots,
   \end{aligned}
  \right. 
\end{equation}
where this solution in $x$ coincides with \eqref{eq:sol1}.

For the system \eqref{eq:3}, K. Okamota constructed its space of initial conditions. His idea is very important (see \cite{O3}). However, singular analysis is very complicated (cf. \cite{MMT,T1}) in the case of Painlev\'e I system \eqref{eq:3}. For example, by his holomorphy condition, the system \eqref{eq:3} can not transformed into a polynomial Hamiltonian system, and is transformed into a complicated rational form (cf. \cite{takano}). In \cite{iwasaki}, K. Iwasaki and S. Okada gave some Hamiltonian structures for the first Painlev\'e system \eqref{eq:3}. They solved this problem by using the algebraic transformation \eqref{eq:5} and its holomorphy condition \eqref{eq:30} (\cite{1}, \cite{Shi}, \cite{iwasaki}).

In this note, we remark that we can do analysis of its accessible singular points for the system \eqref{eq:6}.

Let us consider the regular vector field $V$ associated with the system \eqref{eq:6} defined on $(v,u,t) \in {\mathbb C}^2 \times B$
\begin{equation*}\label{eq:8}
\boxed{%
V=\frac{\partial}{\partial t}+\frac{\partial K}{\partial u}\frac{\partial}{\partial v}-\frac{\partial K}{\partial v}\frac{\partial}{\partial u}
}%
\end{equation*}
to a rational vector field $\tilde V$ on ${\mathbb P}^2 \times B$, where $t \in B={\mathbb C}$.

This rational vector field $\tilde V$ belongs to
\begin{equation*}\label{eq:9}
\tilde V \in H^0({\mathbb P}^2 \times B,\Theta_{{\mathbb P}^2 \times B}(-\log{(H_{{\mathbb P}^2} \times B)})(6(H_{{\mathbb P}^2} \times B))),
\end{equation*}
where the symbol $H_{{\mathbb P}^2} \cong {\mathbb P}^1$ denotes the canonical divisor of ${\mathbb P}^2$ whose self-intersection number of $H_{{\mathbb P}^2}$ is given by $(H_{{\mathbb P}^2})^2=1$.

Since its order of pole is 6, its singularity analysis is difficult. So, we will replace its compactification ${\mathbb P}^2 \times B$ by the Hirzebruch surface of degree four ${\Sigma_4}$ given in next section (see \eqref{eq:10}, Figures 1 and 2). After replacing it, we will see that its order of pole is 1.

After we review the notion of accessible singularity and local index in Section 3, for the system  \eqref{eq:6} we will calculate its accessible singularity and local index in Section 4. In Section 5, we will show that the system \eqref{eq:6} passes the Painlev\'e $\alpha$-test for all accessible singular points $P_i \ (i=1,2,3)$ (see \eqref{eq:17}).

We remark that the system \eqref{eq:6} has two birational symmetries (cf. \cite{iwasaki}, P7 \cite{Shi}, \cite{iwasaki}):
\begin{align}\label{eq:32}
\begin{split}
        &s_{0}:(v,u,t) \rightarrow (-a v,a^4 u,-at),\\
        &s_1: (v,u,t) \rightarrow \left(av,-a^4 \left(u-\frac{t}{v^2}-\frac{4}{v^6} \right),-a t \right),
        \end{split}
        \end{align}
where $a \in \{-1,(-1)^{\frac{1}{5}},-(-1)^{\frac{2}{5}},(-1)^{\frac{3}{5}},-(-1)^{\frac{4}{5}} \}=\{a \in {\mathbb C}|a^5+1=0 \}$. In particular, the restriction ${s_1}|_{a=-1}$ is an automorphism of order 2 for the system \eqref{eq:6} 
\begin{equation}\label{specBa}
{s_1}|_{a=-1}: (v,u,t) \rightarrow \left(-v,- \left(u-\frac{t}{v^2}-\frac{4}{v^6} \right), t \right),
\end{equation}
where $({s_1}|_{a=-1})^2=1$. In \cite{iwasaki}, the above transformations denote $\sigma$.

\begin{figure}
\unitlength 0.1in
\begin{picture}( 52.5900, 40.0700)(  6.1000,-42.8700)
%
\special{pn 8}%
\special{pa 1050 686}%
\special{pa 2636 686}%
\special{fp}%
%
\special{pn 8}%
\special{pa 1064 1572}%
\special{pa 2628 1572}%
\special{fp}%
%
\special{pn 8}%
\special{pa 1336 556}%
\special{pa 1336 1710}%
\special{fp}%
%
\special{pn 8}%
\special{pa 2412 556}%
\special{pa 2412 1704}%
\special{fp}%
\put(27.3700,-12.0200){\makebox(0,0)[lb]{${\Sigma_4}$}}%
%
\special{pn 20}%
\special{pa 1336 686}%
\special{pa 1540 686}%
\special{fp}%
\special{sh 1}%
\special{pa 1540 686}%
\special{pa 1474 666}%
\special{pa 1488 686}%
\special{pa 1474 706}%
\special{pa 1540 686}%
\special{fp}%
%
\special{pn 20}%
\special{pa 1336 686}%
\special{pa 1336 890}%
\special{fp}%
\special{sh 1}%
\special{pa 1336 890}%
\special{pa 1356 822}%
\special{pa 1336 836}%
\special{pa 1316 822}%
\special{pa 1336 890}%
\special{fp}%
\put(13.6300,-6.5000){\makebox(0,0)[lb]{$v$}}%
\put(11.5100,-8.6000){\makebox(0,0)[lb]{$u$}}%
\put(27.3700,-16.1600){\makebox(0,0)[lb]{$H \cong {\mathbb P}^1$}}%
%
\special{pn 20}%
\special{pa 1336 1572}%
\special{pa 1336 1376}%
\special{fp}%
\special{sh 1}%
\special{pa 1336 1376}%
\special{pa 1316 1444}%
\special{pa 1336 1430}%
\special{pa 1356 1444}%
\special{pa 1336 1376}%
\special{fp}%
%
\special{pn 20}%
\special{pa 1330 1572}%
\special{pa 1520 1572}%
\special{fp}%
\special{sh 1}%
\special{pa 1520 1572}%
\special{pa 1452 1552}%
\special{pa 1466 1572}%
\special{pa 1452 1592}%
\special{pa 1520 1572}%
\special{fp}%
%
\special{pn 20}%
\special{pa 2412 686}%
\special{pa 2412 890}%
\special{fp}%
\special{sh 1}%
\special{pa 2412 890}%
\special{pa 2432 822}%
\special{pa 2412 836}%
\special{pa 2392 822}%
\special{pa 2412 890}%
\special{fp}%
%
\special{pn 20}%
\special{pa 2404 686}%
\special{pa 2220 686}%
\special{fp}%
\special{sh 1}%
\special{pa 2220 686}%
\special{pa 2286 706}%
\special{pa 2272 686}%
\special{pa 2286 666}%
\special{pa 2220 686}%
\special{fp}%
%
\special{pn 20}%
\special{pa 2404 1572}%
\special{pa 2404 1368}%
\special{fp}%
\special{sh 1}%
\special{pa 2404 1368}%
\special{pa 2384 1436}%
\special{pa 2404 1422}%
\special{pa 2424 1436}%
\special{pa 2404 1368}%
\special{fp}%
%
\special{pn 20}%
\special{pa 2398 1566}%
\special{pa 2228 1566}%
\special{fp}%
\special{sh 1}%
\special{pa 2228 1566}%
\special{pa 2294 1586}%
\special{pa 2280 1566}%
\special{pa 2294 1546}%
\special{pa 2228 1566}%
\special{fp}%
\put(21.6500,-6.5000){\makebox(0,0)[lb]{$z_1$}}%
\put(24.7100,-8.5300){\makebox(0,0)[lb]{$w_1$}}%
\put(11.1800,-15.1400){\makebox(0,0)[lb]{$w_2$}}%
\put(13.7000,-17.1000){\makebox(0,0)[lb]{$z_2$}}%
\put(24.7100,-15.1400){\makebox(0,0)[lb]{$w_3$}}%
\put(21.9200,-17.1000){\makebox(0,0)[lb]{$z_3$}}%
%
\special{pn 20}%
\special{sh 0.600}%
\special{ar 2412 992 12 12  0.0000000 6.2831853}%
%
\special{pn 20}%
\special{sh 0.600}%
\special{ar 2404 1208 12 14  0.0000000 6.2831853}%
%
\special{pn 20}%
\special{sh 0.600}%
\special{ar 1322 1580 10 12  0.0000000 6.2831853}%
\put(24.7100,-10.5700){\makebox(0,0)[lb]{$P_1$}}%
\put(24.7100,-12.8200){\makebox(0,0)[lb]{$P_2$}}%
\put(11.7200,-17.2500){\makebox(0,0)[lb]{$P_3$}}%
%
\special{pn 8}%
\special{pa 1196 2428}%
\special{pa 2780 2428}%
\special{dt 0.045}%
%
\special{pn 8}%
\special{pa 1482 2298}%
\special{pa 1482 3452}%
\special{dt 0.045}%
\put(28.8300,-29.4400){\makebox(0,0)[lb]{$\tilde{\Sigma}_4$}}%
\put(15.0800,-23.9200){\makebox(0,0)[lb]{$v$}}%
\put(12.9700,-26.0200){\makebox(0,0)[lb]{$u$}}%
%
\special{pn 20}%
\special{sh 0.600}%
\special{ar 2556 2734 12 12  0.0000000 6.2831853}%
%
\special{pn 20}%
\special{sh 0.600}%
\special{ar 2550 2950 12 14  0.0000000 6.2831853}%
\put(26.1600,-27.9900){\makebox(0,0)[lb]{$P_1$}}%
\put(26.1600,-30.2400){\makebox(0,0)[lb]{$P_2$}}%
%
\special{pn 8}%
\special{pa 1956 2364}%
\special{pa 1956 2066}%
\special{fp}%
\special{sh 1}%
\special{pa 1956 2066}%
\special{pa 1936 2132}%
\special{pa 1956 2118}%
\special{pa 1976 2132}%
\special{pa 1956 2066}%
\special{fp}%
%
\special{pn 8}%
\special{pa 1310 3104}%
\special{pa 1746 3292}%
\special{fp}%
\special{pa 1700 3192}%
\special{pa 1550 3598}%
\special{fp}%
\special{pa 1556 3410}%
\special{pa 1760 3714}%
\special{fp}%
\special{pa 1726 3546}%
\special{pa 1468 3836}%
\special{fp}%
\special{pa 1414 3698}%
\special{pa 1706 3882}%
\special{fp}%
\special{pa 1678 3800}%
\special{pa 1258 4280}%
\special{fp}%
%
\special{pn 8}%
\special{pa 2550 2298}%
\special{pa 2556 4288}%
\special{dt 0.045}%
%
\special{pn 8}%
\special{pa 1176 4194}%
\special{pa 2936 4194}%
\special{dt 0.045}%
%
\special{pn 8}%
\special{pa 1550 4062}%
\special{pa 1108 3830}%
\special{fp}%
\special{pa 1176 3946}%
\special{pa 1114 3488}%
\special{fp}%
\special{pa 1188 3604}%
\special{pa 816 3452}%
\special{fp}%
\special{pa 930 3598}%
\special{pa 924 3054}%
\special{fp}%
\special{pa 998 3250}%
\special{pa 610 3060}%
\special{fp}%
%
\special{pn 8}%
\special{pa 670 3190}%
\special{pa 670 2640}%
\special{dt 0.045}%
%
\special{pn 20}%
\special{pa 3754 4156}%
\special{pa 4250 3336}%
\special{fp}%
\special{pa 3882 3292}%
\special{pa 4238 3562}%
\special{fp}%
\special{pa 3876 3402}%
\special{pa 4088 3068}%
\special{fp}%
\special{pa 4094 3168}%
\special{pa 3822 2996}%
\special{fp}%
\special{pa 3842 3068}%
\special{pa 4020 2784}%
\special{fp}%
\special{pa 4054 2872}%
\special{pa 3686 2662}%
\special{fp}%
%
\special{pn 20}%
\special{pa 3910 4070}%
\special{pa 3598 3822}%
\special{fp}%
\special{pa 3714 3742}%
\special{pa 3496 4076}%
\special{fp}%
\special{pa 3570 4040}%
\special{pa 3354 3888}%
\special{fp}%
\special{pa 3442 3808}%
\special{pa 3250 4112}%
\special{fp}%
\special{pa 3338 4106}%
\special{pa 3066 3924}%
\special{fp}%
%
\special{pn 8}%
\special{pa 3208 3344}%
\special{pa 4510 4112}%
\special{dt 0.045}%
%
\special{pn 20}%
\special{pa 4388 4134}%
\special{pa 4354 4124}%
\special{pa 4324 4112}%
\special{pa 4300 4094}%
\special{pa 4282 4068}%
\special{pa 4272 4038}%
\special{pa 4268 4004}%
\special{pa 4272 3970}%
\special{pa 4284 3940}%
\special{pa 4298 3912}%
\special{pa 4320 3886}%
\special{pa 4344 3862}%
\special{pa 4372 3844}%
\special{pa 4402 3832}%
\special{pa 4434 3828}%
\special{pa 4468 3836}%
\special{pa 4496 3856}%
\special{pa 4512 3882}%
\special{pa 4514 3912}%
\special{pa 4510 3940}%
\special{sp}%
%
\special{pn 20}%
\special{pa 4502 3946}%
\special{pa 4482 4040}%
\special{fp}%
\special{sh 1}%
\special{pa 4482 4040}%
\special{pa 4516 3980}%
\special{pa 4494 3988}%
\special{pa 4476 3972}%
\special{pa 4482 4040}%
\special{fp}%
\put(42.5600,-38.0000){\makebox(0,0)[lb]{$\pi, \ (\pi)^2=1$}}%
%
\special{pn 20}%
\special{pa 5374 1776}%
\special{pa 5870 954}%
\special{fp}%
\special{pa 5502 912}%
\special{pa 5856 1180}%
\special{fp}%
\special{pa 5496 1020}%
\special{pa 5706 686}%
\special{fp}%
\special{pa 5714 788}%
\special{pa 5440 616}%
\special{fp}%
\special{pa 5462 686}%
\special{pa 5638 404}%
\special{fp}%
\special{pa 5672 492}%
\special{pa 5306 280}%
\special{fp}%
%
\special{pn 20}%
\special{pa 5496 1754}%
\special{pa 5154 1500}%
\special{fp}%
\special{pa 5250 1442}%
\special{pa 5020 1768}%
\special{fp}%
%
\special{pn 20}%
\special{pa 3856 3858}%
\special{pa 4094 4040}%
\special{fp}%
%
\special{pn 20}%
\special{pa 5468 1434}%
\special{pa 5706 1616}%
\special{fp}%
\put(45.7600,-20.1400){\makebox(0,0)[lb]{Dynkin diagram of type $E_8^{(1)}$}}%
%
\special{pn 8}%
\special{pa 5010 1130}%
\special{pa 4610 1130}%
\special{fp}%
\special{sh 1}%
\special{pa 4610 1130}%
\special{pa 4678 1150}%
\special{pa 4664 1130}%
\special{pa 4678 1110}%
\special{pa 4610 1130}%
\special{fp}%
\put(41.3000,-11.8000){\makebox(0,0)[lb]{$\Sigma_{\epsilon}^{(2)}$}}%
%
\special{pn 8}%
\special{pa 3270 1120}%
\special{pa 4060 1120}%
\special{fp}%
\special{sh 1}%
\special{pa 4060 1120}%
\special{pa 3994 1100}%
\special{pa 4008 1120}%
\special{pa 3994 1140}%
\special{pa 4060 1120}%
\special{fp}%
\put(32.5000,-10.5000){\makebox(0,0)[lb]{Double}}%
\put(32.6000,-13.6000){\makebox(0,0)[lb]{covering}}%
\put(46.0000,-9.8000){\makebox(0,0)[lb]{Blow up at}}%
\put(46.1000,-13.8000){\makebox(0,0)[lb]{eight points}}%
\put(9.9000,-21.4000){\makebox(0,0)[lb]{Blow up at}}%
\put(20.6000,-21.3000){\makebox(0,0)[lb]{twelve points}}%
%
\special{pn 8}%
\special{pa 3050 1120}%
\special{pa 3240 1120}%
\special{dt 0.045}%
\end{picture}%
\label{fig:PIGRSfig4}
\caption{Resolution of accessible singular points and double covering (cf. \cite{iwasaki})}
\end{figure}
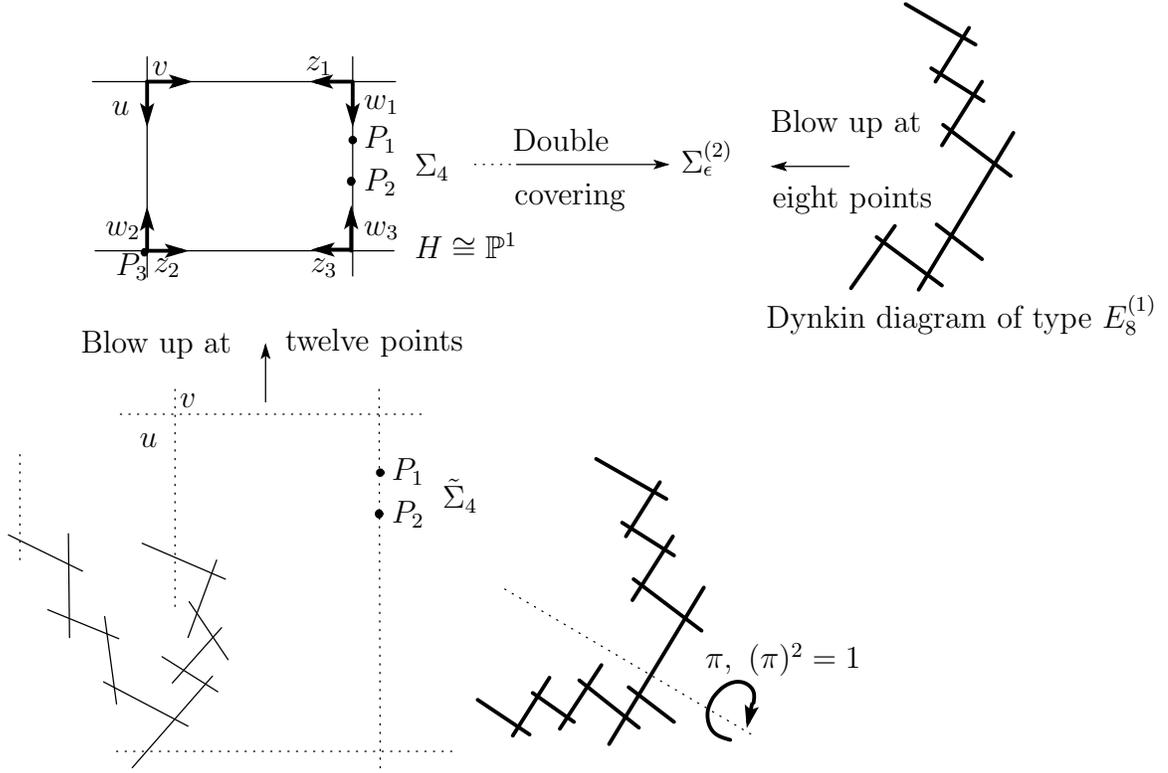

Here, let us consider the holomorphy conditions for the system \eqref{eq:6}.

By resolving the accessible singular point $P_3$ given in Lemma \ref{lem1}, we can obtain the following holomorphy condition for the system \eqref{eq:6}.

We see that the system \eqref{eq:6} becomes again a polynomial Hamiltonian system in the coordinate system $r_3$:(cf. \cite{iwasaki})
\begin{align}\label{eq:30}
r_3:(x_3,y_3)=&\left(v,u-\frac{t}{v^2}-\frac{4}{v^6} \right).
\end{align}
The transformation $r_3$ is birational and symplectic (cf. \cite{iwasaki}, P7 \cite{Shi}).

We note that the condition $r_3$ should be read that $r_3^{-1} \left(K-\frac{1}{v} \right)$ is a polynomial with respect to $x_3,y_3$. In this case, we can obtain
\begin{align}\label{eq:31}
\begin{split}
r_3^{-1} \left(K-\frac{1}{v} \right)=-\frac{x_3^6 y_3^2}{4} + \frac{x_3^5 y_3}{4} - \frac{1}{4} t x_3^4 y_3 + \frac{t x_3^3}{8} - \frac{x_3^4}{16} - 
 \frac{1}{16} t^2 x_3^2 - y_3.
\end{split}
\end{align}

By using this holomorphy condition, we can recover the system  \eqref{eq:6} in a regular vector field $V$:
\begin{equation*}
V=\frac{\partial}{\partial t}+\frac{\partial F}{\partial u}\frac{\partial}{\partial v}-\frac{\partial F}{\partial v}\frac{\partial}{\partial u}, \quad F \in {\mathbb C}(t)[v,u].
\end{equation*}

Let us consider a regular vector field $V$
\begin{equation}
V=\frac{\partial}{\partial t}+\frac{\partial F}{\partial u}\frac{\partial}{\partial v}-\frac{\partial F}{\partial v}\frac{\partial}{\partial u}
\end{equation}
associated with polynomial Hamiltonian system with Hamiltonian $F \in {\mathbb C}(t)[v,u]$. We assume that

$(C1)$ $V \in H^0({\Sigma_4} \times B,\Theta_{{\Sigma_4} \times B}(-\log{((H \cup L) \times B)})((H \cup L) \times B)), \quad B \cong {\mathbb C}$.

$(C2)$ This system becomes again a polynomial Hamiltonian system in the coordinate system $r${\rm : \rm}
\begin{align}\label{holo2}
\begin{split}
r:(X,Y)=&\left(v,u-\frac{t}{v^2}-\frac{4}{v^6} \right).
\end{split}
\end{align}
Then such a system coincides with the Hamiltonian system \eqref{eq:6} with the polynomial Hamiltonian \eqref{eq:7}. Here, the symbol ${\Sigma_4}$ denotes the Hirzebruch surface of degree four ${\Sigma_4}$ given in next section (see \eqref{eq:10}).

We remark that Professor Paul Painlev\'e (see below * in P 346 ;\cite{17},\cite{1,2}) gave its holomorphy condition (cf. \cite{iwasaki}) of the first Painlev\'e system \eqref{eq:3},\eqref{eq:4};
\begin{align}\label{holoP}
\begin{split}
R:(x,y)=&\left(\frac{1}{X^2},\frac{4+t X^4-X^5+2X^6 Y}{2 X^3} \right).
\end{split}
\end{align}
Here, we note that $dy \wedge dx=2dX \wedge dY-d \left(\frac{1}{X} \right) \wedge dt$, where the transformation $R$ is an algebraic transformation of degree 2 (cf. \cite{iwasaki}).

This holomorphy condition \eqref{holoP} can be obtained by composing two transformations \eqref{eq:5}, \eqref{holo2}.

We remark that the holomorphy condition $R$ should be read that
\begin{align*}
\begin{split}
&R \left(H +\frac{1}{X} \right)
\end{split}
\end{align*}
is a polynomial with respect to $X,Y$.

\section{Compactification}

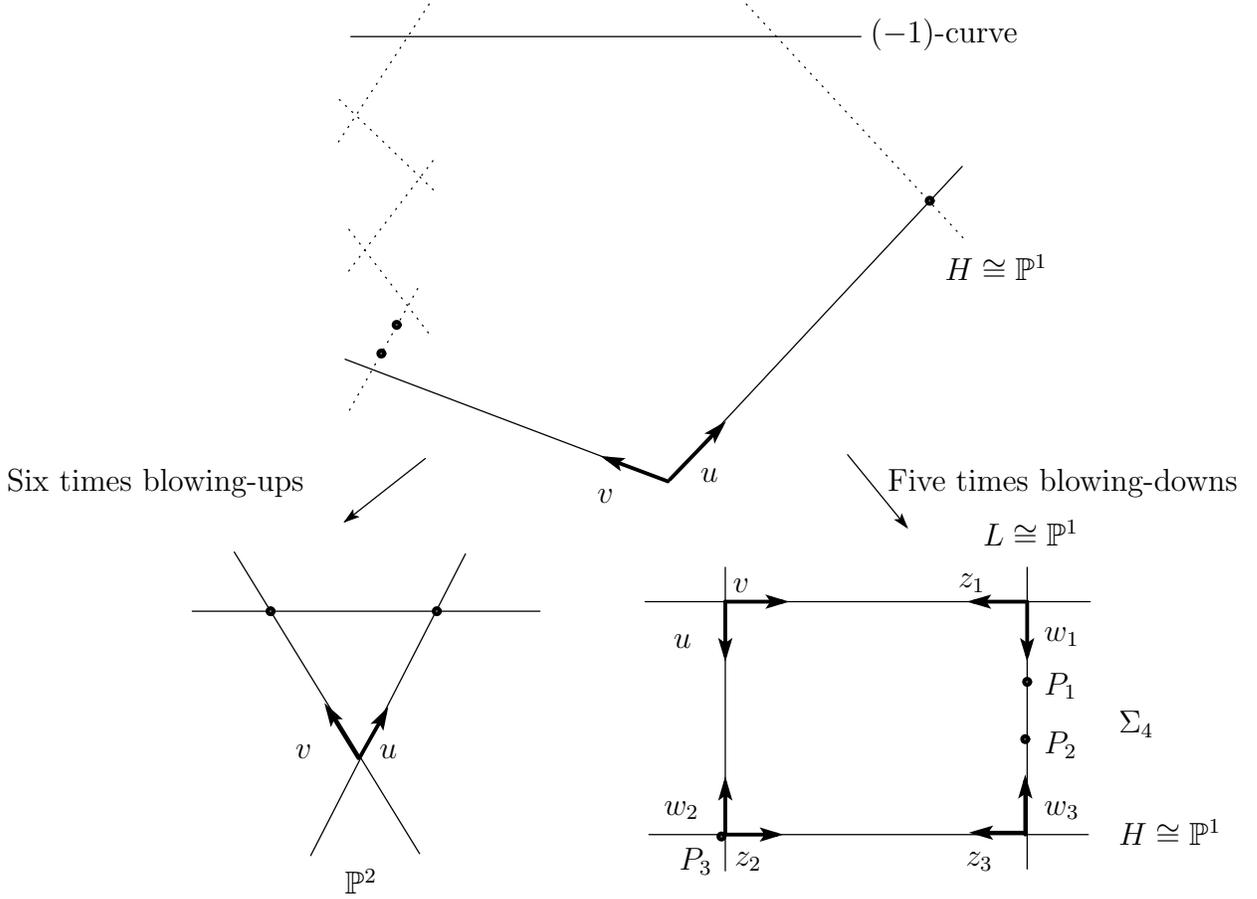
\begin{figure}
\unitlength 0.1in
\begin{picture}( 58.2000, 45.4000)( 12.0000,-47.6000)
%
\special{pn 8}%
\special{pa 4540 3350}%
\special{pa 6870 3350}%
\special{fp}%
%
\special{pn 8}%
\special{pa 4560 4570}%
\special{pa 6860 4570}%
\special{fp}%
%
\special{pn 8}%
\special{pa 4960 3170}%
\special{pa 4960 4760}%
\special{fp}%
%
\special{pn 8}%
\special{pa 6540 3170}%
\special{pa 6540 4750}%
\special{fp}%
\put(70.2000,-40.6000){\makebox(0,0)[lb]{${\Sigma}_4$}}%
%
\special{pn 20}%
\special{pa 4960 3350}%
\special{pa 5260 3350}%
\special{fp}%
\special{sh 1}%
\special{pa 5260 3350}%
\special{pa 5194 3330}%
\special{pa 5208 3350}%
\special{pa 5194 3370}%
\special{pa 5260 3350}%
\special{fp}%
%
\special{pn 20}%
\special{pa 4960 3350}%
\special{pa 4960 3630}%
\special{fp}%
\special{sh 1}%
\special{pa 4960 3630}%
\special{pa 4980 3564}%
\special{pa 4960 3578}%
\special{pa 4940 3564}%
\special{pa 4960 3630}%
\special{fp}%
\put(50.0000,-33.0000){\makebox(0,0)[lb]{$v$}}%
\put(46.9000,-35.9000){\makebox(0,0)[lb]{$u$}}%
\put(70.2000,-46.3000){\makebox(0,0)[lb]{$H \cong {\mathbb P}^1$}}%
%
\special{pn 20}%
\special{pa 4960 4570}%
\special{pa 4960 4300}%
\special{fp}%
\special{sh 1}%
\special{pa 4960 4300}%
\special{pa 4940 4368}%
\special{pa 4960 4354}%
\special{pa 4980 4368}%
\special{pa 4960 4300}%
\special{fp}%
%
\special{pn 20}%
\special{pa 4950 4570}%
\special{pa 5230 4570}%
\special{fp}%
\special{sh 1}%
\special{pa 5230 4570}%
\special{pa 5164 4550}%
\special{pa 5178 4570}%
\special{pa 5164 4590}%
\special{pa 5230 4570}%
\special{fp}%
%
\special{pn 20}%
\special{pa 6540 3350}%
\special{pa 6540 3630}%
\special{fp}%
\special{sh 1}%
\special{pa 6540 3630}%
\special{pa 6560 3564}%
\special{pa 6540 3578}%
\special{pa 6520 3564}%
\special{pa 6540 3630}%
\special{fp}%
%
\special{pn 20}%
\special{pa 6530 3350}%
\special{pa 6260 3350}%
\special{fp}%
\special{sh 1}%
\special{pa 6260 3350}%
\special{pa 6328 3370}%
\special{pa 6314 3350}%
\special{pa 6328 3330}%
\special{pa 6260 3350}%
\special{fp}%
%
\special{pn 20}%
\special{pa 6530 4570}%
\special{pa 6530 4290}%
\special{fp}%
\special{sh 1}%
\special{pa 6530 4290}%
\special{pa 6510 4358}%
\special{pa 6530 4344}%
\special{pa 6550 4358}%
\special{pa 6530 4290}%
\special{fp}%
%
\special{pn 20}%
\special{pa 6520 4560}%
\special{pa 6270 4560}%
\special{fp}%
\special{sh 1}%
\special{pa 6270 4560}%
\special{pa 6338 4580}%
\special{pa 6324 4560}%
\special{pa 6338 4540}%
\special{pa 6270 4560}%
\special{fp}%
\put(61.8000,-33.0000){\makebox(0,0)[lb]{$z_1$}}%
\put(66.3000,-35.8000){\makebox(0,0)[lb]{$w_1$}}%
\put(46.4000,-44.9000){\makebox(0,0)[lb]{$w_2$}}%
\put(50.1000,-47.6000){\makebox(0,0)[lb]{$z_2$}}%
\put(66.3000,-44.9000){\makebox(0,0)[lb]{$w_3$}}%
\put(62.2000,-47.6000){\makebox(0,0)[lb]{$z_3$}}%
%
\special{pn 8}%
\special{pa 2170 3400}%
\special{pa 3990 3400}%
\special{fp}%
\special{pa 2390 3090}%
\special{pa 3360 4670}%
\special{fp}%
\special{pa 3600 3100}%
\special{pa 2790 4680}%
\special{fp}%
%
\special{pn 20}%
\special{pa 3050 4160}%
\special{pa 3180 3930}%
\special{fp}%
\special{sh 1}%
\special{pa 3180 3930}%
\special{pa 3130 3978}%
\special{pa 3154 3976}%
\special{pa 3166 3998}%
\special{pa 3180 3930}%
\special{fp}%
%
\special{pn 20}%
\special{pa 3040 4170}%
\special{pa 2880 3910}%
\special{fp}%
\special{sh 1}%
\special{pa 2880 3910}%
\special{pa 2898 3978}%
\special{pa 2908 3956}%
\special{pa 2932 3956}%
\special{pa 2880 3910}%
\special{fp}%
\put(27.1000,-41.7000){\makebox(0,0)[lb]{$v$}}%
\put(31.5000,-41.7000){\makebox(0,0)[lb]{$u$}}%
%
\special{pn 8}%
\special{pa 3000 390}%
\special{pa 5670 390}%
\special{fp}%
%
\special{pn 8}%
\special{pa 3410 220}%
\special{pa 2930 950}%
\special{dt 0.045}%
\special{pa 2940 720}%
\special{pa 3450 1210}%
\special{dt 0.045}%
\special{pa 3430 1040}%
\special{pa 2980 1630}%
\special{dt 0.045}%
\special{pa 2990 1400}%
\special{pa 3420 1960}%
\special{dt 0.045}%
%
\special{pn 8}%
\special{pa 5070 220}%
\special{pa 6210 1450}%
\special{dt 0.045}%
%
\special{pn 8}%
\special{pa 3350 1710}%
\special{pa 2980 2370}%
\special{dt 0.045}%
%
\special{pn 8}%
\special{pa 2970 2080}%
\special{pa 4660 2720}%
\special{fp}%
\special{pa 4660 2720}%
\special{pa 6200 1070}%
\special{fp}%
%
\special{pn 20}%
\special{pa 4670 2710}%
\special{pa 4940 2430}%
\special{fp}%
\special{sh 1}%
\special{pa 4940 2430}%
\special{pa 4880 2464}%
\special{pa 4904 2468}%
\special{pa 4908 2492}%
\special{pa 4940 2430}%
\special{fp}%
%
\special{pn 20}%
\special{pa 4660 2720}%
\special{pa 4340 2600}%
\special{fp}%
\special{sh 1}%
\special{pa 4340 2600}%
\special{pa 4396 2642}%
\special{pa 4390 2620}%
\special{pa 4410 2606}%
\special{pa 4340 2600}%
\special{fp}%
\put(42.9000,-28.3000){\makebox(0,0)[lb]{$v$}}%
\put(48.3000,-27.2000){\makebox(0,0)[lb]{$u$}}%
%
\special{pn 20}%
\special{sh 0.600}%
\special{ar 3160 2050 16 16  0.0000000 6.2831853}%
%
\special{pn 20}%
\special{sh 0.600}%
\special{ar 3240 1900 16 16  0.0000000 6.2831853}%
%
\special{pn 20}%
\special{sh 0.600}%
\special{ar 6030 1250 16 16  0.0000000 6.2831853}%
%
\special{pn 8}%
\special{pa 3390 2600}%
\special{pa 2970 2930}%
\special{fp}%
\special{sh 1}%
\special{pa 2970 2930}%
\special{pa 3036 2906}%
\special{pa 3012 2898}%
\special{pa 3010 2874}%
\special{pa 2970 2930}%
\special{fp}%
%
\special{pn 8}%
\special{pa 5600 2580}%
\special{pa 5910 2960}%
\special{fp}%
\special{sh 1}%
\special{pa 5910 2960}%
\special{pa 5884 2896}%
\special{pa 5876 2920}%
\special{pa 5852 2922}%
\special{pa 5910 2960}%
\special{fp}%
\put(29.7000,-48.8000){\makebox(0,0)[lb]{${\mathbb P}^2$}}%
%
\special{pn 20}%
\special{sh 0.600}%
\special{ar 6540 3770 16 16  0.0000000 6.2831853}%
%
\special{pn 20}%
\special{sh 0.600}%
\special{ar 6530 4070 16 16  0.0000000 6.2831853}%
%
\special{pn 20}%
\special{sh 0.600}%
\special{ar 4940 4580 16 16  0.0000000 6.2831853}%
\put(57.2000,-4.6000){\makebox(0,0)[lb]{$(-1)$-curve}}%
\put(12.0000,-28.0000){\makebox(0,0)[lb]{Six times blowing-ups}}%
\put(58.1000,-28.0000){\makebox(0,0)[lb]{Five times blowing-downs}}%
\put(66.3000,-38.6000){\makebox(0,0)[lb]{$P_1$}}%
\put(66.3000,-41.7000){\makebox(0,0)[lb]{$P_2$}}%
\put(47.2000,-47.8000){\makebox(0,0)[lb]{$P_3$}}%
%
\special{pn 20}%
\special{sh 0.600}%
\special{ar 2580 3400 16 16  0.0000000 6.2831853}%
%
\special{pn 20}%
\special{sh 0.600}%
\special{ar 3450 3400 16 16  0.0000000 6.2831853}%
\put(61.1000,-16.7000){\makebox(0,0)[lb]{$H \cong {\mathbb P}^1$}}%
\put(63.1000,-30.6000){\makebox(0,0)[lb]{$L \cong {\mathbb P}^1$}}%
\end{picture}%
\label{fig:PIGRSfig1}
\caption{Each symbol $\bullet$ denotes accessible singular point.}
\end{figure}

In order to consider the singularity analysis for the system \eqref{eq:6}, as a compactification of ${\mathbb C}^2$ which is the phase space of the system \eqref{eq:6}, we take the following Hirzebruch surface of degree four ${\Sigma_4}$, which is obtained by gluing four copies of ${\mathbb C}^2$ via the following identification (see Figures 2 and 3):
\begin{align}\label{eq:10}
\begin{split}
&{\Sigma_4}=U_0 \cup \bigcup_{i=1}^{3} U_j, \quad U_j \cong {\mathbb C}^2 \ni (z_j,w_j) \ (j=0,1,2,3),\\
&z_0=v, \ w_0=u, \quad z_1=\frac{1}{v}, \ w_1=\left(\left(uv-\frac{1}{2} \right)v-\frac{t}{2} \right)v^2,\\
&z_2=z_0, \ w_2=\frac{1}{w_0}, \quad z_3=z_1, \ w_3=\frac{1}{w_1}.
\end{split}
\end{align}
We define two divisors $H$ and $L$ (see Figure 2) on ${\Sigma_4}$:
\begin{align}\label{eq:11}
\begin{split}
&H:=\{(z_2,w_2) \in U_2|w_2=0\} \cup \{(z_3,w_3) \in U_3|w_3=0\} \cong {\mathbb P}^1,\\
&L:=\{(z_1,w_1) \in U_1|z_1=0\} \cup \{(z_3,w_3) \in U_3|z_3=0\} \cong {\mathbb P}^1.
\end{split}
\end{align}
Each self-intersection number of $H$ and $L$ is given by
\begin{equation}\label{eq:12}
(H)^2=4, \quad (L)^2=0.
\end{equation}
After a series of successive six times blowing-ups and five times blowing-downs on projective surface ${\mathbb P}^2$ {\rm  (see Figure 2), \rm} we obtain Hirzebruch surface of degree four ${\Sigma_4}$ and a birational morphism $\varphi:{\Sigma_4} \cdots \rightarrow {\mathbb P}^2$. Its canonical divisor $K_{\Sigma_4}$ is given by
\begin{equation}\label{eq:13}
K_{\Sigma_4}=-2H,
\end{equation}
where the symbol $H$ denotes the proper transform of $H$ by $\varphi$.

On ${\Sigma_4} \times B$ in \eqref{eq:10}, we see that this rational vector field $\tilde V$ associated with the system \eqref{eq:6} belongs to
\begin{equation}
\boxed{%
\tilde V \in H^0({\Sigma_4} \times B,\Theta_{{\Sigma_4} \times B}(-\log{((H \cup L) \times B)})((H \cup L) \times B)).
}%
\end{equation}
We remark that this rational vector field $\tilde V$ has a pole along the divisors $H$ and $L$, whose order is one.

\section{Review of accessible singularity and local index}
Let us review the notion of {\it accessible singularity}. Let $B$ be a connected open domain in $\mathbb C$ and $\pi : {\mathcal W} \longrightarrow B$ a smooth proper holomorphic map. We assume that ${\mathcal H} \subset {\mathcal W}$ is a normal crossing divisor which is flat over $B$. Let us consider a rational vector field $\tilde v$ on $\mathcal W$ satisfying the condition
\begin{equation*}
\tilde v \in H^0({\mathcal W},\Theta_{\mathcal W}(-\log{\mathcal H})({\mathcal H})).
\end{equation*}
Fixing $t_0 \in B$ and $P \in {\mathcal W}_{t_0}$, we can take a local coordinate system $(x_1,\ldots ,x_n)$ of ${\mathcal W}_{t_0}$ centered at $P$ such that ${\mathcal H}_{\rm smooth \rm}$ can be defined by the local equation $x_1=0$.
Since $\tilde v \in H^0({\mathcal W},\Theta_{\mathcal W}(-\log{\mathcal H})({\mathcal H}))$, we can write down the vector field $\tilde v$ near $P=(0,\ldots ,0,t_0)$ as follows:
\begin{equation*}
\tilde v= \frac{\partial}{\partial t}+g_1 
\frac{\partial}{\partial x_1}+\frac{g_2}{x_1} 
\frac{\partial}{\partial x_2}+\cdots+\frac{g_n}{x_1} 
\frac{\partial}{\partial x_n}.
\end{equation*}
This vector field defines the following system of differential equations
\begin{equation}\label{39}
\frac{dx_1}{dt}=g_1(x_1,\ldots,x_n,t),\ \frac{dx_2}{dt}=\frac{g_2(x_1,\ldots,x_n,t)}{x_1},\cdots, \frac{dx_n}{dt}=\frac{g_n(x_1,\ldots,x_n,t)}{x_1}.
\end{equation}
Here $g_i(x_1,\ldots,x_n,t), \ i=1,2,\ldots ,n,$ are holomorphic functions defined near $P$.

\begin{definition}\label{Def1}
With the above notation, assume that the rational vector field $\tilde v$ on $\mathcal W$ satisfies the condition
$$
(A) \quad \tilde v \in H^0({\mathcal W},\Theta_{\mathcal W}(-\log{\mathcal H})({\mathcal H})).
$$
We say that $\tilde v$ has an {\it accessible singularity} at $P=(0,\dots ,0,t_0)$ if
\begin{equation}
\boxed{%
x_1=0 \ {\rm and \rm} \ g_i(0,\ldots,0,t_0)=0 \ {\rm for \rm} \ {\rm every \rm} \ i, \ 2 \leq i \leq n.
}%
\end{equation}
\end{definition}

If $P \in {\mathcal H}_{{\rm smooth \rm}}$ is not an accessible singularity, all solutions of the ordinary differential equation passing through $P$ are vertical solutions, that is, the solutions are contained in the fiber ${\mathcal W}_{t_0}$ over $t=t_0$. If $P \in {\mathcal H}_{\rm smooth \rm}$ is an accessible singularity, there may be a solution of \eqref{39} which passes through $P$ and goes into the interior ${\mathcal W}-{\mathcal H}$ of ${\mathcal W}$.

Here we review the notion of {\it local index}. Let $v$ be an algebraic vector field with an accessible singular point $\overrightarrow{p}=(0,\ldots,0)$ and $(x_1,\ldots,x_n)$ be a coordinate system in a neighborhood centered at $\overrightarrow{p}$. Assume that the system associated with $v$ near $\overrightarrow{p}$ can be written as

{\Small
\begin{align}\label{b}
\begin{split}
\frac{d}{dt}\begin{pmatrix}
             x_1 \\
             x_2 \\
             \vdots\\
             x_{n-1} \\
             x_n
             \end{pmatrix}=\frac{1}{x_1}\left\{\begin{bmatrix}
             a_{11} & 0 & 0 & \hdots & 0 \\
             a_{21} & a_{22} & 0 &  \hdots & 0 \\
             \vdots & \vdots & \ddots & 0 & 0 \\
             a_{(n-1)1} & a_{(n-1)2} & \hdots & a_{(n-1)(n-1)} & 0 \\
             a_{n1} & a_{n2} & \hdots & a_{n(n-1)} & a_{nn}
             \end{bmatrix}\begin{pmatrix}
             x_1 \\
             x_2 \\
             \vdots\\
             x_{n-1} \\
             x_n
             \end{pmatrix}+\begin{pmatrix}
             x_1h_1(x_1,\ldots,x_n,t) \\
             h_2(x_1,\ldots,x_n,t) \\
             \vdots\\
             h_{n-1}(x_1,\ldots,x_n,t) \\
             h_n(x_1,\ldots,x_n,t)
             \end{pmatrix}\right\},\\
              (h_i \in {\mathbb C}(t)[x_1,\ldots,x_n], \ a_{ij} \in {\mathbb C}(t))
             \end{split}
             \end{align}}
where $h_1$ is a polynomial which vanishes at $\overrightarrow{p}$ and $h_i$, $i=2,3,\ldots,n$ are polynomials of order at least 2 in $x_1,x_2,\ldots,x_n$, We call ordered set of the eigenvalues $(a_{11},a_{22},\cdots,a_{nn})$ {\it local index} at $\overrightarrow{p}$.

We are interested in the case with local index
\begin{equation}\label{integer}
\left(1,\frac{a_{22}(t)}{a_{11}(t)},\ldots,\frac{a_{nn}(t)}{a_{11}(t)} \right) \in {\mathbb Z}^{n}.
\end{equation}

If each component of $\left(1,\frac{a_{22}(t)}{a_{11}(t)},\ldots,\frac{a_{nn}(t)}{a_{11}(t)} \right)$ has the same sign, we may resolve the accessible singularity by blowing-up finitely many times. However, when different signs appear, we may need to both blow up and blow down.
\begin{center}
\begin{tabular}{|c||c|c|} \hline 
& $\left(\frac{a_{22}(t)}{a_{11}(t)},\ldots,\frac{a_{nn}(t)}{a_{11}(t)} \right)$  & Resolution of accessible sing. \\ \hline
Positive sign & ${\mathbb N}^{n-1}$ & Blowing-up \\ \hline
Different signs &  ${\mathbb Z}^{n-1}$ & both Blow up and Blow down  \\ \hline
\end{tabular}
\end{center}

The $\alpha$-test,
\begin{equation}\label{poiuy}
t=t_0+\alpha T, \quad x_i=\alpha X_i, \quad \alpha \rightarrow 0,
\end{equation}
yields the following reduced system:
\begin{align}\label{ppppppp}
\begin{split}
\frac{d}{dT}\begin{pmatrix}
             X_1 \\
             X_2 \\
             \vdots\\
             X_{n-1} \\
             X_n
             \end{pmatrix}=\frac{1}{X_1}\begin{bmatrix}
             a_{11}(t_0) & 0 & 0 & \hdots & 0 \\
             a_{21}(t_0) & a_{22}(t_0) & 0 &  \hdots & 0 \\
             \vdots & \vdots & \ddots & 0 & 0 \\
             a_{(n-1)1}(t_0) & a_{(n-1)2}(t_0) & \hdots & a_{(n-1)(n-1)}(t_0) & 0 \\
             a_{n1}(t_0) & a_{n2}(t_0) & \hdots & a_{n(n-1)}(t_0) & a_{nn}(t_0)
             \end{bmatrix}\begin{pmatrix}
             X_1 \\
             X_2 \\
             \vdots\\
             X_{n-1} \\
             X_n
             \end{pmatrix},
             \end{split}
             \end{align}
where $a_{ij}(t_0) \in {\mathbb C}$. Fixing $t=t_0$, this system is the system of the first order ordinary differential equation with constant coefficient. Let us solve this system. At first, we solve the first equation:
\begin{equation}
X_1(T)=a_{11}(t_0)T+C_1 \quad (C_1 \in {\mathbb C}).
\end{equation}
Substituting this into the second equation in \eqref{ppppppp}, we can obtain the first order linear ordinary differential equation:
\begin{equation}
\frac{dX_2}{dT}=\frac{a_{22}(t_0) X_2}{a_{11}(t_0)T+C_1}+a_{21}(t_0).
\end{equation}
By variation of constant, in the case of $a_{11}(t_0) \not= a_{22}(t_0)$ we can solve explicitly:
\begin{equation}
X_2(T)=C_2(a_{11}(t_0)T+C_1)^{\frac{a_{22}(t_0)}{a_{11}(t_0)}}+\frac{a_{21}(t_0)(a_{11}(t_0)T+C_1)}{a_{11}(t_0)-a_{22}(t_0)} \quad (C_2 \in {\mathbb C}).
\end{equation}
This solution is a single-valued solution if and only if
$$
\frac{a_{22}(t_0)}{a_{11}(t_0)} \in {\mathbb Z}-\{1\}.
$$
In the case of $a_{11}(t_0)=a_{22}(t_0)$ we can solve explicitly:
\begin{equation}
X_2(T)=C_2(a_{11}(t_0)T+C_1)+\frac{a_{21}(t_0)(a_{11}(t_0)T+C_1){\rm Log}(a_{11}(t_0)T+C_1)}{a_{11}(t_0)} \quad (C_2 \in {\mathbb C}).
\end{equation}
This solution is a single-valued solution if and only if
$$
a_{21}(t_0)=0.
$$
Of course, $\frac{a_{22}(t_0)}{a_{11}(t_0)}=1 \in {\mathbb Z}$.
In the same way, we can obtain the solutions for each variables $(X_3,\ldots,X_n)$.

\begin{tabular}{|c|} \hline 
The conditions $\frac{a_{jj}(t)}{a_{11}(t)} \in {\mathbb Z}, \ (j=2,3,\ldots,n)$ are necessary condition in order to have\\
the Painlev\'e property. \\ \hline 
\end{tabular}

\begin{center}
\begin{tabular}{|c||c|c|} \hline 
& $\left(\frac{a_{22}(t)}{a_{11}(t)},\ldots,\frac{a_{nn}(t)}{a_{11}(t)} \right)$  & Movable singularities \\ \hline
Painlev\'e type & ${\mathbb Z}$ & Only pole  \\ \hline
Other Non-Linear Equation &  ${\mathbb Q},{\mathbb R}$ and ${\mathbb C}$ & Algebraic sing. or others  \\ \hline
\end{tabular}
\end{center}

For example, we consider the Painlev\'e VI equation. The sixth Painlev\'e equation is equivalent to the following Hamiltonian system:
\begin{equation}\label{PVI}
  \left\{
  \begin{aligned}
   \frac{dx}{dt} =&\frac{\partial H_{VI}}{\partial y}=\frac{1}{t(t-1)}\{2y(x-t)(x-1)x-(\alpha_0-1)(x-1)x\\
&-\alpha_3(x-t)x-\alpha_4(x-t)(x-1)\},\\
   \frac{dy}{dt} =&-\frac{\partial H_{VI}}{\partial x}=\frac{1}{t(t-1)}[-\{(x-t)(x-1)+(x-t)x+(x-1)x\}y^2\\
&+\{(\alpha_0-1)(2x-1)+\alpha_3(2x-t)+\alpha_4(2x-t-1)\}y\\
&-\alpha_2(\alpha_1+\alpha_2)]
   \end{aligned}
  \right. 
\end{equation}
with the polynomial Hamiltonian
\begin{align}\label{HVI}
\begin{split}
&H_{VI}(x,y,t;\alpha_0,\alpha_1,\alpha_2,\alpha_3,\alpha_4)\\
&=\frac{1}{t(t-1)}[y^2(x-t)(x-1)x-\{(\alpha_0-1)(x-1)x+\alpha_3(x-t)x\\
&+\alpha_4(x-t)(x-1)\}y+\alpha_2(\alpha_1+\alpha_2)x]  \quad (\alpha_0+\alpha_1+2\alpha_2+\alpha_3+\alpha_4=1).
\end{split}
\end{align}
Since each right hand side of this system is polynomial with respect to $x,y$, by Cauchy's existence and uniqueness theorem of solutions, there exists unique holomorphic solution with initial values $(x,y)=(x_0,y_0) \in {\mathbb C}^2$.

Let us extend the regular vector field defined on ${\mathbb C}^2 \times B$
$$
v=\frac{\partial}{\partial t}+\frac{\partial H_{VI}}{\partial y}\frac{\partial}{\partial x}-\frac{\partial H_{VI}}{\partial x}\frac{\partial}{\partial y}
$$
to a rational vector field on ${\Sigma_{-\alpha_2}^{(2)}} \times B$, where $B={\mathbb C}-\{0,1\}$.

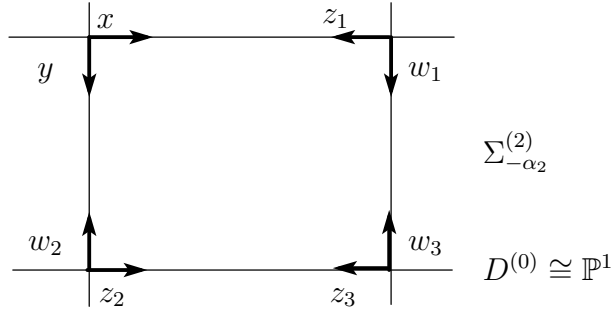
\begin{figure}
\unitlength 0.1in
\begin{picture}( 24.8000, 16.3000)( 19.9000,-20.0000)
%
\special{pn 8}%
\special{pa 1990 590}%
\special{pa 4320 590}%
\special{fp}%
%
\special{pn 8}%
\special{pa 2010 1810}%
\special{pa 4310 1810}%
\special{fp}%
%
\special{pn 8}%
\special{pa 2410 410}%
\special{pa 2410 2000}%
\special{fp}%
%
\special{pn 8}%
\special{pa 3990 410}%
\special{pa 3990 1990}%
\special{fp}%
\put(44.7000,-13.0000){\makebox(0,0)[lb]{${\Sigma_{-\alpha_2}^{(2)}}$}}%
%
\special{pn 20}%
\special{pa 2410 590}%
\special{pa 2710 590}%
\special{fp}%
\special{sh 1}%
\special{pa 2710 590}%
\special{pa 2644 570}%
\special{pa 2658 590}%
\special{pa 2644 610}%
\special{pa 2710 590}%
\special{fp}%
%
\special{pn 20}%
\special{pa 2410 590}%
\special{pa 2410 870}%
\special{fp}%
\special{sh 1}%
\special{pa 2410 870}%
\special{pa 2430 804}%
\special{pa 2410 818}%
\special{pa 2390 804}%
\special{pa 2410 870}%
\special{fp}%
\put(24.5000,-5.4000){\makebox(0,0)[lb]{$x$}}%
\put(21.4000,-8.3000){\makebox(0,0)[lb]{$y$}}%
\put(44.7000,-18.7000){\makebox(0,0)[lb]{$D^{(0)} \cong {\mathbb P}^1$}}%
%
\special{pn 20}%
\special{pa 2410 1810}%
\special{pa 2410 1540}%
\special{fp}%
\special{sh 1}%
\special{pa 2410 1540}%
\special{pa 2390 1608}%
\special{pa 2410 1594}%
\special{pa 2430 1608}%
\special{pa 2410 1540}%
\special{fp}%
%
\special{pn 20}%
\special{pa 2400 1810}%
\special{pa 2680 1810}%
\special{fp}%
\special{sh 1}%
\special{pa 2680 1810}%
\special{pa 2614 1790}%
\special{pa 2628 1810}%
\special{pa 2614 1830}%
\special{pa 2680 1810}%
\special{fp}%
%
\special{pn 20}%
\special{pa 3990 590}%
\special{pa 3990 870}%
\special{fp}%
\special{sh 1}%
\special{pa 3990 870}%
\special{pa 4010 804}%
\special{pa 3990 818}%
\special{pa 3970 804}%
\special{pa 3990 870}%
\special{fp}%
%
\special{pn 20}%
\special{pa 3980 590}%
\special{pa 3710 590}%
\special{fp}%
\special{sh 1}%
\special{pa 3710 590}%
\special{pa 3778 610}%
\special{pa 3764 590}%
\special{pa 3778 570}%
\special{pa 3710 590}%
\special{fp}%
%
\special{pn 20}%
\special{pa 3980 1810}%
\special{pa 3980 1530}%
\special{fp}%
\special{sh 1}%
\special{pa 3980 1530}%
\special{pa 3960 1598}%
\special{pa 3980 1584}%
\special{pa 4000 1598}%
\special{pa 3980 1530}%
\special{fp}%
%
\special{pn 20}%
\special{pa 3970 1800}%
\special{pa 3720 1800}%
\special{fp}%
\special{sh 1}%
\special{pa 3720 1800}%
\special{pa 3788 1820}%
\special{pa 3774 1800}%
\special{pa 3788 1780}%
\special{pa 3720 1800}%
\special{fp}%
\put(36.3000,-5.4000){\makebox(0,0)[lb]{$z_1$}}%
\put(40.8000,-8.2000){\makebox(0,0)[lb]{$w_1$}}%
\put(20.9000,-17.3000){\makebox(0,0)[lb]{$w_2$}}%
\put(24.6000,-20.0000){\makebox(0,0)[lb]{$z_2$}}%
\put(40.8000,-17.3000){\makebox(0,0)[lb]{$w_3$}}%
\put(36.7000,-20.0000){\makebox(0,0)[lb]{$z_3$}}%
\end{picture}%
\label{fig:E6figure1}
\caption{Rational surface ${\Sigma_{-\alpha_2}^{(2)}}$}
\end{figure}

Here, we review the rational surface ${\Sigma_{-\alpha_2}^{(2)}}$, which is obtained by gluing four copies of ${\mathbb C}^2$ via the following identification:
\begin{align}
\begin{split}
&U_j \cong {\mathbb C}^2 \ni (z_j,w_j) \ (j=0,1,2,3),\\
&z_0=x, \ w_0=y, \quad z_1=\frac{1}{x}, \ w_1=-(xy+\alpha_2)x,\\
&z_2=z_0, \ w_2=\frac{1}{w_0}, \quad z_3=z_1, \ w_3=\frac{1}{w_1}.
\end{split}
\end{align}
We define a divisor $D^{(0)}$ on ${\Sigma_{-\alpha_2}^{(2)}}$:
\begin{equation}
D^{(0)}=\{(z_2,w_2) \in U_2|w_2=0\} \cup \{(z_3,w_3) \in U_3|w_3=0\} \cong {\mathbb P}^1.
\end{equation}
The self-intersection number of $D^{(0)}$ is given by
\begin{equation*}
(D^{(0)})^2=2.
\end{equation*}
In the coordinate system $(z_1,w_1)$ the right hand side of this system is polynomial with respect to $z_1,w_1$. However, on the boundary divisor $D^{(0)} \cong {\mathbb P}^1$ this system has a pole in each coordinate system $(z_i,w_i) \ i=2,3$. By calculating the accessible singular points on $D^{(0)}$, we obtain simple four singular points $z_2=0,1,t,\infty$ (see Definition \ref{Def1}).

By rewriting the system at each singular point, this rational vector field has a pole along the divisor $D^{(0)}$, whose order is one.

By resolving all singular points, we can construct the space of initial conditions of the Painlev\'e VI system. This space parametrizes all meromorphic solutions including holomorphic solutions.

 Conversely, we can recover the Painlev\'e VI system by all patching data of its space of initial conditions. At first, we decompose its patching data into the pair of singular points and local index around each singular point.

Now, let us rewrite the system centered at each singular point $X=0,1,t,\infty$.

1. By taking the coordinate system $(X,Y)=(z_2,w_2)$ centered at the point $(z_2,w_2)=(0,0)$, the system is given by

{\Small
\begin{align*}
&\frac{d}{dt}\begin{pmatrix}
             X \\
             Y 
             \end{pmatrix}=\frac{1}{t(t-1)Y} \{ t\begin{pmatrix}
             2 & -\alpha_4   \\
             0 & 1
             \end{pmatrix}\begin{pmatrix}
             X \\
             Y 
             \end{pmatrix}+\begin{pmatrix}
             -2(t+1) & \alpha_0-1+\alpha_4+t(\alpha_3+\alpha_4) \\
             0 & -2(t+1) 
             \end{pmatrix}\begin{pmatrix}
             X^2 \\
             XY 
             \end{pmatrix}\\
&+\begin{pmatrix}
             2 & -(\alpha_0+\alpha_3+\alpha_4)+1 \\
             0 & 3 
             \end{pmatrix}\begin{pmatrix}
             X^3 \\
             X^2 Y 
             \end{pmatrix} \}+\begin{pmatrix}
             0 \\
             \frac{-\{(\alpha_0-1)(2X-1)+\alpha_3(2X-t)+\alpha_4(2X-t-1)\}Y+\alpha_2(\alpha_1+\alpha_2)Y^2}{t(t-1)}
             \end{pmatrix}.
             \end{align*}}
Now, let us make a change of variables $X,Y,t$ with a small parameter $\alpha$:
\begin{equation}
X=\alpha Z, \quad Y=\alpha W, \quad t=t_0+\alpha T \quad (t_0 \in {\mathbb C}-\{0,1\}).
\end{equation}
Then the system can also be written in the new variables $Z,W,T$. This new system tends to the system as $\alpha \rightarrow 0$
\begin{align}\label{ZZZZ}
\frac{d}{dT}\begin{pmatrix}
             Z \\
             W 
             \end{pmatrix}&=\frac{1}{W}\left\{\begin{pmatrix}
             \frac{2}{t_0-1} & -\frac{\alpha_4}{t_0-1}  \\
             0 & \frac{1}{t_0-1} 
             \end{pmatrix}\begin{pmatrix}
             Z \\
             W 
             \end{pmatrix}\right\}.
             \end{align}
Fixing $t=t_0$, this system is the system of the first order ordinary differential equation with constant coefficient. Let us solve this system. At first, we solve the second equation:
\begin{equation}
W(T)=\frac{T}{t_0-1}+C_1 \quad (C_1 \in {\mathbb C}).
\end{equation}
Substituting this into the first equation in \eqref{ZZZZ}, we can obtain the first order linear ordinary differential equation:
\begin{equation}
\frac{dZ}{dT}=\frac{t_0-1}{T+C_1(t_0-1)}\left(\frac{2}{t_0-1}Z-\frac{\alpha_4}{t_0-1} \left(\frac{T}{t_0-1}+C_1 \right) \right).
\end{equation}
By variation of constant, we can solve explicitly:
\begin{equation}
Z(T)=C_2\{T+(t_0-1)C_1 \}^2+\frac{\alpha_4(T+(t_0-1)C_1)}{t_0-1} \quad (C_2 \in {\mathbb C}).
\end{equation}
Thus, we can obtain single-valued solutions. For the Painlev\'e property, this is the necessary condition.

\vspace{0.5cm}
In the same way, we can obtain the following:

2. By taking the coordinate system $(X,Y)=(z_2-1,w_2)$ centered at the point $(z_2,w_2)=(1,0)$, the system is given by
\begin{align*}
\frac{d}{dt}\begin{pmatrix}
             X \\
             Y 
             \end{pmatrix}&=\frac{1}{Y}\left\{\begin{pmatrix}
             -\frac{2}{t} & \frac{\alpha_3}{t}  \\
             0 & -\frac{1}{t} 
             \end{pmatrix}\begin{pmatrix}
             X \\
             Y 
             \end{pmatrix}+\cdots\right\}
             \end{align*}

3. By taking the coordinate system $(X,Y)=(z_2-t,w_2)$ centered at the point $(z_2,w_2)=(t,0)$, the system is given by
\begin{align*}
\frac{d}{dt}\begin{pmatrix}
             X \\
             Y 
             \end{pmatrix}&=\frac{1}{Y}\left\{\begin{pmatrix}
             2 & -\alpha_0  \\
             0 & 1 
             \end{pmatrix}\begin{pmatrix}
             X \\
             Y 
             \end{pmatrix}+\cdots\right\}
             \end{align*}

4. By taking the coordinate system $(X,Y)=(z_3,w_3)$ centered at the point $(z_3,w_3)=(0,0)$, the system is given by
\begin{align*}
\frac{d}{dt}\begin{pmatrix}
             X \\
             Y 
             \end{pmatrix}&=\frac{1}{Y}\left\{\begin{pmatrix}
             \frac{2}{t(t-1)} & -\frac{\alpha_1}{t(t-1)}  \\
             0 & \frac{1}{t(t-1)} 
             \end{pmatrix}\begin{pmatrix}
             X \\
             Y 
             \end{pmatrix}+\cdots\right\}.
             \end{align*}
Thus, we have proved that the Hamiltonian system \eqref{PVI},\eqref{HVI} passes the Painlev\'e $\alpha$-test for all accessible singular points $X=0,1,t,\infty$

\begin{equation*}
\begin{pmatrix}
X=0 & X=1 & X=t & X=\infty\\
\frac{1}{t-1}\begin{pmatrix}
2 & -\alpha_4\\
0 & 1 
\end{pmatrix} & -\frac{1}{t}\begin{pmatrix}
2 & -\alpha_3\\
0 & 1 
\end{pmatrix} & \begin{pmatrix}
2 & -\alpha_0\\
0 & 1 
\end{pmatrix} & \frac{1}{t(t-1)}\begin{pmatrix}
2 & -\alpha_1\\
0 & 1 
\end{pmatrix}
\end{pmatrix}.
\end{equation*}

\begin{center}
\begin{tabular}{|c|} \hline 
The pair of accessible singular points and matrix of linear approximation\\
around each point is called {\it Painlev\'e scheme}. \\ \hline 
\end{tabular}
\end{center}

\section{Accessible singularities and Local index for our system}

For the system \eqref{eq:6}, let us calculate its accessible singularities.

Around the point $(z_1,w_1)=(0,0)$, the system can be rewritten as follows:
\begin{equation}\label{eq:15}
  \left\{
  \begin{aligned}
   \frac{dz_1}{dt} &=-z_1^2 + \frac{w_1}{2},\\
   \frac{dw_1}{dt} &=\frac{2 t - w_1^2}{2 z_1} + 4 z_1 w_1,
   \end{aligned}
  \right. 
\end{equation}
and around the point $(z_2,w_2)=(0,0)$, the system can be rewritten as follows:
\begin{equation}\label{eq:16}
  \left\{
  \begin{aligned}
   \frac{dz_2}{dt} &=1 + \frac{t z_2^4}{4} + \frac{z_2^5}{4} - \frac{z_2^6}{2 w_2},\\
   \frac{dw_2}{dt} &=-\frac{3 z_2^5}{2} + t z_2^3 w_2 + \frac{5 z_2^4 w_2}{4} - \frac{1}{8} t^2 z_2 w_2^2 - 
 \frac{3}{8} t z_2^2 w_2^2 - \frac{z_2^3 w_2^2}{4},
   \end{aligned}
  \right. 
\end{equation}
and around the point $(z_3,w_3)=(0,0)$, the system can be rewritten as follows:
\begin{equation}\label{eq:z3eq}
  \left\{
  \begin{aligned}
   \frac{dz_3}{dt} &=-z_3^2 + \frac{1}{2w_3},\\
   \frac{dw_3}{dt} &=\frac{1 - 2t w_3^2}{2 z_3} - 4 z_3 w_3.
   \end{aligned}
  \right. 
\end{equation}

\begin{figure}[h]
\unitlength 0.1in
\begin{picture}( 24.8000, 16.3000)( 14.9000,-22.5000)
%
\special{pn 8}%
\special{pa 1490 840}%
\special{pa 3820 840}%
\special{fp}%
%
\special{pn 8}%
\special{pa 1510 2060}%
\special{pa 3810 2060}%
\special{fp}%
%
\special{pn 8}%
\special{pa 1910 660}%
\special{pa 1910 2250}%
\special{fp}%
%
\special{pn 8}%
\special{pa 3490 660}%
\special{pa 3490 2240}%
\special{fp}%
\put(39.7000,-15.5000){\makebox(0,0)[lb]{${\Sigma}_4$}}%
%
\special{pn 20}%
\special{pa 1910 840}%
\special{pa 2210 840}%
\special{fp}%
\special{sh 1}%
\special{pa 2210 840}%
\special{pa 2144 820}%
\special{pa 2158 840}%
\special{pa 2144 860}%
\special{pa 2210 840}%
\special{fp}%
%
\special{pn 20}%
\special{pa 1910 840}%
\special{pa 1910 1120}%
\special{fp}%
\special{sh 1}%
\special{pa 1910 1120}%
\special{pa 1930 1054}%
\special{pa 1910 1068}%
\special{pa 1890 1054}%
\special{pa 1910 1120}%
\special{fp}%
\put(19.5000,-7.9000){\makebox(0,0)[lb]{$v$}}%
\put(16.4000,-10.8000){\makebox(0,0)[lb]{$u$}}%
\put(39.7000,-21.2000){\makebox(0,0)[lb]{$H \cong {\mathbb P}^1$}}%
%
\special{pn 20}%
\special{pa 1910 2060}%
\special{pa 1910 1790}%
\special{fp}%
\special{sh 1}%
\special{pa 1910 1790}%
\special{pa 1890 1858}%
\special{pa 1910 1844}%
\special{pa 1930 1858}%
\special{pa 1910 1790}%
\special{fp}%
%
\special{pn 20}%
\special{pa 1900 2060}%
\special{pa 2180 2060}%
\special{fp}%
\special{sh 1}%
\special{pa 2180 2060}%
\special{pa 2114 2040}%
\special{pa 2128 2060}%
\special{pa 2114 2080}%
\special{pa 2180 2060}%
\special{fp}%
%
\special{pn 20}%
\special{pa 3490 840}%
\special{pa 3490 1120}%
\special{fp}%
\special{sh 1}%
\special{pa 3490 1120}%
\special{pa 3510 1054}%
\special{pa 3490 1068}%
\special{pa 3470 1054}%
\special{pa 3490 1120}%
\special{fp}%
%
\special{pn 20}%
\special{pa 3480 840}%
\special{pa 3210 840}%
\special{fp}%
\special{sh 1}%
\special{pa 3210 840}%
\special{pa 3278 860}%
\special{pa 3264 840}%
\special{pa 3278 820}%
\special{pa 3210 840}%
\special{fp}%
%
\special{pn 20}%
\special{pa 3480 2060}%
\special{pa 3480 1780}%
\special{fp}%
\special{sh 1}%
\special{pa 3480 1780}%
\special{pa 3460 1848}%
\special{pa 3480 1834}%
\special{pa 3500 1848}%
\special{pa 3480 1780}%
\special{fp}%
%
\special{pn 20}%
\special{pa 3470 2050}%
\special{pa 3220 2050}%
\special{fp}%
\special{sh 1}%
\special{pa 3220 2050}%
\special{pa 3288 2070}%
\special{pa 3274 2050}%
\special{pa 3288 2030}%
\special{pa 3220 2050}%
\special{fp}%
\put(31.3000,-7.9000){\makebox(0,0)[lb]{$z_1$}}%
\put(35.8000,-10.7000){\makebox(0,0)[lb]{$w_1$}}%
\put(15.9000,-19.8000){\makebox(0,0)[lb]{$w_2$}}%
\put(19.6000,-22.5000){\makebox(0,0)[lb]{$z_2$}}%
\put(35.8000,-19.8000){\makebox(0,0)[lb]{$w_3$}}%
\put(31.7000,-22.5000){\makebox(0,0)[lb]{$z_3$}}%
%
\special{pn 20}%
\special{sh 0.600}%
\special{ar 3490 1260 16 16  0.0000000 6.2831853}%
%
\special{pn 20}%
\special{sh 0.600}%
\special{ar 3480 1560 16 16  0.0000000 6.2831853}%
%
\special{pn 20}%
\special{sh 0.600}%
\special{ar 1890 2070 16 16  0.0000000 6.2831853}%
\put(35.8000,-13.5000){\makebox(0,0)[lb]{$P_1$}}%
\put(35.8000,-16.6000){\makebox(0,0)[lb]{$P_2$}}%
\put(16.7000,-22.7000){\makebox(0,0)[lb]{$P_3$}}%
\end{picture}%
\label{fig:PIGRSfig2}
\caption{Each symbol $\bullet$ denotes accessible singular point.}
\end{figure}
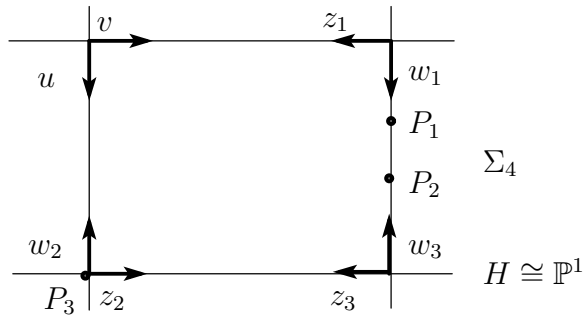

\begin{lemma}\label{lem1}
The rational vector field $\tilde V$ associated with the system \eqref{eq:6}  has three accessible singular points $P_i \ (i=1,2,3)$ \rm{(see figure 4)}$:$
\begin{equation}\label{eq:17}
  \left\{
  \begin{aligned}
    P_1=&\{(z_1,w_1)|z_1=0, \ w_1=\sqrt{2t}\},\\
    P_2=&\{(z_1,w_1)|z_1=0, \ w_1=-\sqrt{2t}\},\\
    P_3=&\{(z_2,w_2)|z_2=w_2=0\},
   \end{aligned}
  \right. 
\end{equation}
where the point $P_3$ has multiplicity of order 6.
\end{lemma}

This lemma can be proven by a direct calculation. \qed

We see that the system \eqref{eq:15} is invariant under the following birational transformation $\pi$:
\begin{equation}\label{autom}
\pi : (z_1,w_1) \rightarrow (-z_1,4z_1^2-w_1).
\end{equation}
This transformation $\pi$ changes two accessible singular points $P_1$ and $P_2$.

We note that pulling back the transformation ${s_1}|_{a=-1}$ in \eqref{specBa} by the birational transformation $(z_1,w_1)=\left(\frac{1}{v},\left(\left(uv-\frac{1}{2} \right)v-\frac{t}{2} \right)v^2 \right)$, we can obtain the above transformation \eqref{autom}.

Next let us calculate its local index at $P_1$ and $P_2$.
\begin{center}
\begin{tabular}{|c|c|c|c|} \hline 
Singular point & Type of local index & Resonance  \\ \hline 
$P_1$ & $(\frac{\sqrt{2t}}{2},-\sqrt{2t})$ & $\frac{-\sqrt{2t}}{\frac{\sqrt{2t}}{2}}=-2$  \\ \hline 
$P_2$ & $(\frac{\sqrt{2t}}{2},-\sqrt{2t})$ & $\frac{-\sqrt{2t}}{\frac{\sqrt{2t}}{2}}=-2$  \\ \hline 
\end{tabular}
\end{center}
We see that the rational vector field $\tilde V$ associated with the system \eqref{eq:6} has negative resonance $-2$ at each accessible singular point $P_1$ and $P_2$, respectively (see \eqref{negaseries}).

\section{Painlev\'e $\alpha$-test}
In this section, we will show that the system \eqref{eq:6} passes the Painlev\'e $\alpha$-test for all accessible singular points $P_i \ (i=1,2,3)$.

\begin{proposition}\label{th:1}
The rational vector field $\tilde V$ associated with the system \eqref{eq:6} passes the Painlev\'e $\alpha$-test at each of accessible singular points $P_i \ (i=1,2,3)$.
\end{proposition}

{\it Proof.} Around the point $P_1$, the system can be rewritten in the coordinate system $(X,Y)=(z_1,w_1-\sqrt{2t})$:
\begin{equation}\label{eq:18}
  \left\{
  \begin{aligned}
   \frac{dX}{dt} &=\frac{\sqrt{2t}}{2} + \frac{Y}{2}-X^2,\\
   \frac{dY}{dt} &=-\frac{\sqrt{2t} Y }{X}-\frac{Y^2 }{2X} + 4 XY+4\sqrt{2t} X-\frac{1}{\sqrt{2t}}.
   \end{aligned}
  \right. 
\end{equation}

We remark that the relations between the coordinate system $(x,y)$ in \eqref{eq:3} and the coordinate system $(X,Y)$ in \eqref{eq:18} are given as follows:
\begin{equation}
  \left\{
  \begin{aligned}
   x &=X^2,\\
   y &=-X ( 2 X^2 - Y-\sqrt{2t} )
   \end{aligned}
  \right. 
\end{equation}
and
\begin{equation}
  \left\{
  \begin{aligned}
   X &=\sqrt{x},\\
   Y &=2x-\sqrt{2t}+\frac{y}{\sqrt{x}}.
   \end{aligned}
  \right. 
\end{equation}

The $\alpha$-test,
\begin{equation}\label{eq:19}
t=t_0+\alpha T, \quad X=\alpha X_1, \ Y=\alpha Y_1, \quad \alpha \rightarrow 0,
\end{equation}
yields the following reduced system:
\begin{equation}\label{eq:20}
   \frac{dX_1}{dT} =\frac{\sqrt{2t_0}}{2}, \quad \frac{dY_1}{dT} =-\frac{\sqrt{2t_0} Y_1}{X_1}-\frac{1}{\sqrt{2t_0}}, \quad (t_0 \in {\mathbb C}).
\end{equation}
We remark that this system is a system of the first-order ordinary differential equations with {\it constant} coefficients.

Solving this system, we can obtain its solution:
\begin{equation}\label{eq:21}
  \left\{
  \begin{aligned}
  X_1[T] &= \frac{\sqrt{2t_0}}{2}T + C_1,\\
   Y_1[T] &=\frac{-\sqrt{2}t_0 T^3 - 6 C_1 \sqrt{t_0} T^2- 6 \sqrt{2} C_1^2 T+ 
  3C_2 \sqrt{t_0}}{3 \sqrt{t_0} (\sqrt{2t_0}T + 2 C_1)^2} \quad (C_1,C_2 \in {\mathbb C}),
   \end{aligned}
  \right. 
\end{equation}
where $C_1$ and $C_2$ are integral constants.

Since these solutions are rational solutions in $T$, we see that the system \eqref{eq:20} can be solved by single-valued solutions.

At the accessible singular point $P_1$, the system \eqref{eq:6} passes the Painlev\'e $\alpha$-test.

By the same way, we can prove in the case of the accessible singular point $P_2$.

Next, we consider the accessible singular point  $P_3$. This point $P_3$ has multiplicity of order 6:
\begin{align}\label{eq:22}
\begin{split}
\frac{d}{dt}\begin{pmatrix}
             z_2 \\
             w_2 
             \end{pmatrix}=&\frac{1}{w_2} \{ \begin{pmatrix}
             0 & 1  \\
             0 & 0
             \end{pmatrix}\begin{pmatrix}
             z_2 \\
             w_2 
             \end{pmatrix}+\begin{pmatrix}
             0 & 0  \\
             0 & 0
             \end{pmatrix} \begin{pmatrix}
             z_2^2 \\
             z_2w_2 
             \end{pmatrix}+\begin{pmatrix}
             0 & 0  \\
             0 & 0
             \end{pmatrix} \begin{pmatrix}
             z_2^3 \\
             z_2^2 w_2 
             \end{pmatrix}\\
             &+\begin{pmatrix}
             0 & 0  \\
             0 & 0
             \end{pmatrix} \begin{pmatrix}
             z_2^4 \\
             z_2^3 w_2 
             \end{pmatrix}+\begin{pmatrix}
             0 & 0  \\
             0 & 0
             \end{pmatrix}\begin{pmatrix}
             z_2^5 \\
             z_2^4 w_2 
             \end{pmatrix}+\begin{pmatrix}
             -\frac{1}{2} & 0  \\
             0 & -\frac{3}{2}
             \end{pmatrix}\begin{pmatrix}
             z_2^6 \\
             z_2^5 w_2 
             \end{pmatrix}+\cdots \}.
             \end{split}
             \end{align}
The eigenvalues of five matrices in this Painlev\'e expansion around the point $P_3$ are all zero.

By doing successive six times blowing-ups, this accessible singular point transforms into a {\it simple} singular point, and we can get the coordinate system (cf. \cite{MMT}):
\begin{equation}\label{eq:23}
(X,Y):=(v,u v^6).
\end{equation}
In the coordinate system $(X,Y)=(v,u v^6)$, the system can be rewritten as follows:
\begin{equation}\label{eq:24}
  \left\{
  \begin{aligned}
   \frac{dX}{dt} &=1 - \frac{Y}{2} + \frac{t X^4}{4} + \frac{X^5}{4},\\
   \frac{dY}{dt} &=-\frac{3 (Y - 4) Y}{2 X} + \frac{t^2 X^7}{8} + \frac{3 t X^8}{8} + \frac{X^9}{4} + 
 \frac{1}{2} t X^3 Y + \frac{X^4 Y}{4}.
   \end{aligned}
  \right. 
\end{equation}
Taking into account of $Y \not= 0$ when $X=0$ (cf. \cite{MMT}), we see that the system \eqref{eq:24} has the accessible singular point:
\begin{equation}\label{eq:25}
\tilde{P}_3=\{(X_,Y):=(0,4)\}.
\end{equation}
Around the point $\tilde{P}_3$, the system \eqref{eq:24} can be rewritten in the coordinate system $(X_1,Y_1)=(X,Y-4)$:
\begin{equation}\label{eq:26}
  \left\{
  \begin{aligned}
   \frac{dX_1}{dt} &=-1 - \frac{Y_1}{2} + \frac{t X_1^4}{4} + \frac{X_1^5}{4},\\
   \frac{dY_1}{dt} &=-\frac{6Y_1}{X_1}-\frac{3Y_1^2}{2X_1}+ 
  2 t X_1^3  + X_1^4+ \frac{t^2 X_1^7}{8} + \frac{3 t X_1^8}{8} + \frac{X_1^9}{4} + 
 \frac{1}{2} t X_1^3 Y_1 + \frac{X_1^4 Y_1}{4}.
   \end{aligned}
  \right. 
\end{equation}
The $\alpha$-test,
\begin{equation}\label{eq:27}
t=t_0+\alpha T, \quad X_1=\alpha X_2, \ Y_1=\alpha Y_2, \quad \alpha \rightarrow 0,
\end{equation}
yields the following reduced system:
\begin{equation}\label{eq:28}
   \frac{dX_2}{dT} =-1, \quad \frac{dY_2}{dT} =-\frac{6Y_2}{X_2}.
\end{equation}
Solving this system, we can obtain its solution:
\begin{equation}\label{eq:29}
  X_2[T] = -(T- C_1), \quad Y_2[T] =C_2(T - C_1)^6 \quad (C_1,C_2 \in {\mathbb C}),
\end{equation}
where $C_1$ and $C_2$ are integral constants.

Since these solutions are polynomial solutions in $T$, we see that the system \eqref{eq:28} can be solved by single-valued solutions.

At the accessible singular point $P_3$, the system \eqref{eq:6} passes the Painlev\'e $\alpha$-test.

Thus, we have completed the proof of Proposition \ref{th:1}. \qed

We remark that in the system \eqref{eq:26} we can resolve its accessible singular point by six times successive blowing-ups (cf. \cite{Shi,iwasaki}).

\end{document}